\documentclass[leqno,11pt]{amsart}
\usepackage{fourier}
\usepackage[margin=1in]{geometry}
\usepackage{amssymb}
\usepackage{amsfonts}
\usepackage{amsmath}
\usepackage{fancyhdr}
\usepackage{subfig}
\usepackage{graphicx}
\usepackage{caption}
\usepackage{latexsym}
\usepackage{url}
\usepackage{hyperref}
\usepackage[capitalize]{cleveref}
\usepackage{mathrsfs}
\usepackage{color,xcolor}
\usepackage{enumitem}
\usepackage{amsthm}
\usepackage{indentfirst}
\usepackage[sort&compress,square,numbers]{natbib}

\theoremstyle{plain}
\newtheorem{theorem}{Theorem}[section]
\newtheorem{lemma}[theorem]{Lemma}

\newtheorem{proposition}[theorem]{Proposition}
\newtheorem{hyp}{Condition}
\theoremstyle{definition}

\newtheorem{remark}[theorem]{Remark}

\numberwithin{equation}{section}
\allowdisplaybreaks

\newcommand\norm[1]{\left\lVert#1\right\rVert}

\def \Im{\mathop{\rm Im}}
\def \Re{\mathop{\rm Re}}

\newcommand{\rd}{\mathrm{d}}

\newcommand{\eps}{\varepsilon}
\newcommand{\veps}{\varepsilon}

\newcommand{\hess}{\mathrm{Hess}}

\DeclareMathOperator{\tr}{tr}
\numberwithin{equation}{section}

\begin{document}

\title[Asymptotic analysis of diabatic surface hopping algorithm in the adiabatic and non-adiabatic limits]{Asymptotic analysis of diabatic surface hopping algorithm in the adiabatic and non-adiabatic limits} 

\author{Zhenning Cai} \address{ Department of Mathematics, National University of Singapore, Singapore
119076.} \email{matcz@nus.edu.sg}

\author{Di Fang} \address{ Department of Mathematics, Simons Institute for the Theory of Computing, and Challenge Institute for Quantum Computation, University of
  California, Berkeley, Berkeley CA 94720, USA} \email{difang@berkeley.edu}

\author{Jianfeng Lu} \address{Departments of Mathematics, Physics, and
  Chemistry, Duke University, Durham NC 27708, USA}
\email{jianfeng@math.duke.edu} 

\date{\today}

\begin{abstract}
 Surface hopping algorithms, as an important class of quantum dynamics simulation algorithms for non-adiabatic dynamics, are typically performed in the adiabatic representation, which can break down in the presence of ill-defined adiabatic potential energy surfaces (PESs) and adiabatic coupling term. Another issue of surface hopping algorithms is the difficulty in capturing the correct scaling of the transition rate in the Marcus (weak-coupling/non-adiabatic) regime. Though the first issue can be circumvented by exploiting the diabatic representation, diabatic surface hopping algorithms usually lack justification on the theoretical level. We consider the diabatic surface hopping algorithm proposed in [Fang, Lu. Multiscale Model. Simul. 16:4, 1603-1622, 2018] and provide the asymptotic analysis of the transition rate in the Marcus regime that justifies the correct scaling for the spin-boson model. We propose two conditions that guarantee the correctness for general potentials. In the opposite (strong-coupling/adiabatic) regime, we derive the asymptotic behavior of the algorithm that interestingly matches a type of mean-field description. The techniques used here may shed light on the analysis for other diabatic-based algorithms. 
\end{abstract}

\maketitle

\section{Introduction}

The Born-Oppenheimer approximation \cite{Born-oppenheimer}, with the assumption that the electronic degrees of freedom adjust instantaneously to the dynamics of nuclei, has been a basic and key treatment in simplifying and simulating  chemical systems. Fundamental as it is, it considers the motion of nuclei governed by a single adiabatic potential energy surface (PES), which unfortunately becomes inadequate in describing a great number of chemical events where non-adiabatic effects play an crucial role, such as photochemical processes, electron transfer, ultrafast laser experiments, electrochemical reactions and so on. To tackle this issue, different approaches are proposed -- the fully quantum type such as
the multi-configuration time-dependent Hartree \cite{hartree3,hartree1,hartree2}, the semiclassical type such as Meyer-Miller dynamics \cite{miller1,miller2} and multiple spawning methods \cite{spawn}, and the mixed quantum-classical methods such as the Ehrenfest dynamics and surface hopping approaches, to name a few. We refer interested readers to \cite{tully2012perspective,conical_review,ehrenfest_review} for a overview. Among them, surface hopping algorithms, first proposed by John Tully \cite{tully1971, tully1990} are a type of Monte-Carlo algorithms that incorporate the non-adiabatic transitions as a stochastic hop in between different potential energy surfaces, and has become one of the popular tools in simulating quantum systems with non-adiabatic phenomena.

To apply the idea of surface hopping, one starts with a certain choice of representations for the fast degrees of freedom, namely the adiabatic or diabatic representations. While the former fantastically works for a great many problems, it could give rise to singular or large non-adiabatic coupling terms in some cases, which motivates the alternative choice -- the diabatic picture that possesses smooth diabatic coupling terms even near the transition zone. Indeed, it has been pointed out \cite{nitzan2013chemical} that even for the simple Landau-Zener problem, the diabatic surfaces are more appealing in the weak-coupling/high speed limit (\textit{aka} the non-adiabatic limit), while in the large coupling/slow motion limit (\textit{aka} the adiabatic limit) the adiabatic ones are more desirable. Beside the potential numerical advantange, sometimes the diabatic basis may also appear more physical. As is described in \cite{tully2012perspective}, a common example is the spin transitions \cite{spin} where the diabatic PESs are of pure spin character and transitions induced by the spin-orbit interaction while the adiabatic PESs correspond to mixed spin character. Another example is the electonic transfer between two centers witnessed by a nuclear coordinates of a solvent dipole \cite{nitzan2013chemical}. For further discussions on diabatic states and their use, we refer interested readers to a review on this subject \cite{diabatic_annal}. 

The challenge in the diabatic picture, however, is a well-justified algorithm. 
Physical and chemical intuitions can sometimes play a crucial role in the error correction and algorithm designing to capture the correct transition rates, but it would be nice to provide a systematic justification. 
In our previous work \cite{fang-lu}, we proposed a surface hopping algorithm in diabatic representations, based on time dependent perturbation theory and semiclassical analysis of the matrix Schr\"odinger equation. In the spirit of the Tully's fewest switches surface hopping (FSSH), the algorithm can be viewed as a Monte Carlo sampling algorithm on the semiclassical path space for piecewise deterministic path with stochastic jumps between the energy surfaces. The algorithm is proved to produce an accurate wave function as a trajectory sample average in the diabatic representations under mild assumptions on the PESs. Besides the wavefunction, the algorithm is numerically observed to capture the correct transition rates. To be specific, in the weak-coupling (non-adiabatic) regime, the algorithm is numerically observed to capture the correct asymptotic limit for the transition rates that scales quadratically with the off diagonal term, which agrees with the celebrated Marcus formula \cite{marcus1956theory}, while as is pointed out in \cite{landry2011communication} that the original Tully's fewest switch in the adiabatic picture gives linear scaling, which is unphysical. In the opposite regime, however, when the size of the off-diagonal coupling term is very large, a lot more trajectories are needed for the trajectory average to converge in practice. Motivated by these numerical evidence, it is of our interests to justify the observations in a more rigorous manner. In particular, the focus of this paper is to investigate the behavior of the diabatic surface hopping algorithm under these two asymptotic limits. 

These two asymptotic limit is by itself an interesting mathematical regimes as aforementioned. 
From numerical analysis view point, an investigation of the limiting regimes provides a better understanding of the algorithm, and help to identify the difficulty $\delta \to \infty$, where the Monte Carlo algorithm becomes difficult to converge. Making clear of the limiting behavior could help us for improving algorithm design utilizing some multi-scale ideas, such as applying the Heterogeneous Multiscale Method (HMM)~\cite{EEngquist2003,EEngquistLiRen2007,AbdulleEEngquistVanden-Eijnden2012}, averaging algorithms, etc.

When it comes to analysis and tests of non-adiabatic techniques, the spin-boson model is typically considered as the testbed for various algorithms, thanks to its simple form of the harmonic diabatic PESs. The spin-boson model, describing a qubit (system) coupled with a finite-temperature bosonic bath, has the Hamiltonian
\[
\hat H = -\delta \sigma_x + \frac{1}{2} c \sigma_z (b^\dag + b) + w b^\dag b, 
\]
where $\sigma_x$ and $\sigma_z$ are standard Pauli matrices, the $b^\dag$ and $b$ are creation and annihilation operators with frequency $\omega$, and $2\delta$ is hopping rate usually denoted in physical literature. In other words, in the first quantization picture, its diabatic potential matrix is given by harmonic oscillators which take the form 
\begin{equation} \label{eq:V_spin_boson}
V(x) = \frac{1}{2}\omega^2x^2 + c\begin{pmatrix}
x & 0\\ 0& -x
\end{pmatrix} - \delta \begin{pmatrix}
0 & 1 \\ 1 & 0
\end{pmatrix}
\end{equation}
in the one-dimensional case, and the high dimensional case will be addressed later. 
The special form of the harmonic diabatic PESs can help to simplify the analysis and clarify the underlying mechanism, in particular in the semiclassical description the quadratic potential can sometimes diminish the asymptotic errors in the semiclassical ansatz. This unfortunately may lead to misleading performance (see, e.g., \cite{AndoSanter2003,RyabinkinHsiehKapralIzmaylov2014} for the absence of the geometric phase effects). Motivated by the spin-boson model and its possible limitations, we also consider the two-level system in the diabatic representation with general diabatic PESs, where the diabatic potential matrix is 
\begin{equation} \label{eq:V_matrix}
V(x) = 
\begin{pmatrix}
  V_{0}(x) & \delta \\
  \delta & V_{1}(x)
\end{pmatrix},
\end{equation}
where $V_0(x)$ and $V_1(x)$ are assumed to be smooth but not necessarily quadratic.

\vspace{1em}
\noindent\textbf{Contribution:}
The first contribution of this work is to prove analytically that in the non-adiabatic (weak-coupling) limit, that is, the regime of $\delta \ll 1$, the the diabatic surface hopping algorithm proposed in \cite{fang-lu} indeed predicts the correct scaling factor for the Marcus' gold-rule rate, as has been observed numerically. Our proof relies on a careful investigation of the stationary phase analysis. In particular, we propose a set of generic non-degeneracy conditions for the underlying Hessian matrix, under which the transition rate is proved to have the correct scaling as in the Marcus formula. The first condition is a generic condition on the potential energy surfaces, while the other regards the change of the position with respect to the jumping time. Both conditions are physically intuitive, and we verify that the spin-boson model indeed satisfies such non-degeneracy condition.

The second contribution regards the opposite regime when $\delta \to \infty$, which is in fact a non-perturbative regime in the diabatic picture. For the spin-boson model we show the wave function converges to its mean dynamics that contain fast oscillations and behave as an Ehrenfest evolution on an average PES of the two diabatic PESs. In this regime, the perturbation perspective can not be applied. We consider the stochastic representation instead and derive the limiting dynamics using a two-scale expansion for the probability description and the result matches the heuristics. 

\vspace{1em}

\noindent\textbf{Related works:} 
The investigation and application of the fewest-switch-type trajectory surface hopping, starting from Tully's seminal work~\cite{tully1990}, has been a widely studied subject (see \cite{WangAkimovPrezhdo2016,SubotnikJainEtAl2016,conical_review} for reviews of recent progresses and applications), but the mathematical justification has only been provided recently in \cite{lu2016frozen} for the adiabatic representation and \cite{fang-lu} for the diabatic representation. 

The high-level idea of the mathematical framework is to express the part of dynamics that is hard to deal with (the non-adiabatic dynamics in this case) as a perturbation of some easy dynamics (the dynamics governed by the PESs here). Treating coupling terms as the forcing term and applying the Duhamel's principle iteratively, one represents the solution into a Dyson series, which is then interpreted and realized as a stochastic ensemble average of the trajectory information. Such strategy
has been employed in other stochastic algorithms for quantum dynamics simulation, such as path integral thermal dynamics~\cite{LuZhou2017,LuZhou2018,LuLuZhou2020}, quantum kinetic Monte Carlo algorithm~\cite{CaiLu2018}, bold-line Monte Carlo method~\cite{gull2010bold}, continuous-time Monte Carlo method~\cite{RubtsovSavkinLichtenstein2005,GullMillisEtAl2011} and inchworm Monte Carlo method~\cite{cohen2015taming,CaiLuYang2020,CaiLuYang2021,CaiLuYang2022}.

In our case, the trajectory average is based on a generalization of the Herman-Kluk propagator. The transition rate can be viewed as an observable expectation. 
The numerical analysis of the observable expectation has been studied in \cite{LasserSattlegger2017,LasserLubich2020} for the scalar Schr\"odinger equation with a single PES.

We point out that besides the fewest-switch-type surface hopping, there is another branch of trajectory surface hopping algorithms -- the Landau-Zener type -- where each trajectory hopping only once (see, e.g., \cite{Lasser2005,lasser2007,fer_las,BelyaevLasserTrigila2014,avoided_lasser,BelyaevDomckeLasser2015}). It is interesting to observe that the diabatic fewest-switch type surface hopping can be approximately reduced to a single-hop protocol in the weak-coupling regime.

 \vspace{1em}
 
\noindent\textbf{Organization:} 
The rest of this paper is now organized as follows: In \cref{sec:alg_review} we introduces the problem set-up of the spin-boson model and revisit the diabatic surface hopping algorithm proposed in \cite{fang-lu} including its deterministic and stochastic representations. The non-adiabatic limit (weak-coupling) is considered in \cref{sec:weak}. We show the asymptotic scaling using principle of the stationary phase of the deterministic Herman-Kluk ansatz and propose a non-degeneracy condition of the Hessian matrix for general diabatic PESs. 
The verification that the spin-boson model indeed satisfies the non-degeneracy condition is postponed to \cref{sec:appendix}.
\cref{sec:large} is devoted to an investigation of the adiabatic limit (large coupling). The mean dynamics is derived using a two-scale expansion to the Kolmogorov forward equation of the stochastic dynamics. Finally, in \cref{sec:conclusion}, we conclude with some further remarks.

\section{Preliminaries on the diabatic surface hopping algorithm} \label{sec:alg_review}
In this section, we consider a two-level Schr\"odinger equation and introduce the diabatic surface hopping algorithm proposed in our previous work \cite{fang-lu}. 

\subsection{The Schr\"odinger equation in the diabatic representation}
Consider the two-level Schr\"{o}dinger equation (in the diabatic representation):
\begin{equation}
i\veps \partial _{t}%
\begin{pmatrix}
u_{0} \\ 
u_{1}%
\end{pmatrix}%
=-\frac{\veps^2}{2}\Delta _{x}%
\begin{pmatrix}
u_{0} \\ 
u_{1}%
\end{pmatrix}%
+%
\begin{pmatrix}
V_{0} & \delta \\ 
\delta & V_{1}%
\end{pmatrix}%
\begin{pmatrix}
u_{0} \\ 
u_{1}%
\end{pmatrix}%
,  \label{eq:schd}
\end{equation}%
with initial condition 
\begin{equation}
u(0,x)=u_{\text{in}}(x)=%
\begin{pmatrix}
u_{0}(0,x) \\ 
0%
\end{pmatrix}%
,  \label{initial_data}
\end{equation}%
where $\left( t,x\right) \in \mathbb{R}^{+}\times \mathbb{R}^{d}$,
$u\left( t,x\right) =\left[ u_{0}\left( t,x\right), \ u_{1}\left(
    t,x\right) \right]^T$
is the (two component) wave function, and we have assumed without loss
of generality that the initial condition concentrates on the first
energy surface. The Hamiltonian is given by $H=-\frac{%
  \varepsilon ^{2}}{2}\Delta _{x}+V\left( x\right)$, where $V(x)$ is a matrix potential as defined in \eqref{eq:V_matrix}. 
In particular, the spin-boson model as defined in \eqref{eq:V_spin_boson} has the
diabatic potential energy surfaces (PESs) in the form of $d$-dimensional harmonic oscillators $V_0(x) =\frac{1}{2}\omega^2|x|^2 + c\cdot x $ and $V_1(x) =\frac{1}{2}\omega^2 |x|^2 - c \cdot x $. 

It is helpful to define the adiabatic PESs to understand the underlying non-adiabatic dynamics, though we will not use the adiabatic representation in the algorithm. The adiabatic PESs are defined as the eigenvalues of the matrix $V(x)$, namely,
\[
\frac{V_0 + V_1}{2} \pm \sqrt{ \frac{(V_0-V_1)^2}{4} + \delta^2 }.
\]
The gap of adiabatic PESs is given by their difference, which is
\[
2 \sqrt{ \frac{(V_0-V_1)^2}{4} + \delta^2 } =  \sqrt{ {(V_0-V_1)^2} + 4\delta^2 }.
\]
The values of $x$ such that $\lvert V_0(x) -  V_1(x) \rvert$ achieves the minimum are known as the avoided crossings (or hopping points), at which the gap achieves it minimum value. 

In this setup, we have two parameters: $\varepsilon$ is the
semiclassical parameter and $\delta$ is a scaling parameter for the
amplitude of the off-diagonal terms in the matrix potential. One particularly interesting asymptotic limit is the Marcus regime (non-adiabatic limit) when $\delta \to 0$ with fixed $\veps$, which is related to the celebrated Marcus theory \cite{marcus1956theory}. In this regime, the gap of the adiabatic PESs
 is small. Another is the adiabatic limit when the diabatic coupling term $\delta \to \infty$, in which case the adiabatic PESs has large gap
\[
\sqrt{ (V_0 - V_1)^2 + 4\delta^2} \geq \lvert \delta \rvert \to \infty.
\]

\subsection{The diabatic surface hopping algorithm}
We consider a deterministic ansatz to the two-level Schr\"odinger equation, which is a generalization of the Herman-Kluk propagator, aka the frozen Gaussian approximation (FGA), of the scalar Schr\"odinger equation. The deterministic ansatz of total wave function is a summation of the contributions of all number of hops:
\begin{equation}  \label{eq:u_series}
u_\text{F}\left( t,x\right) =%
\begin{pmatrix}
u_\text{F}^{\left( 0\right) }\left( t,x\right) +u_\text{F}^{\left( 2\right) }\left(
t,x\right) +u_\text{F}^{\left( 4\right) }\left( t,x\right) +\cdots  \\ 
u_\text{F}^{\left( 1\right) }\left( t,x\right) +u_\text{F}^{\left( 3\right) }\left(
t,x\right) +u_\text{F}^{\left( 5\right) }\left( t,x\right) +\cdots 
\end{pmatrix}%
,
\end{equation}%
where the contribution initiated at surface $0$, with $n$ hoppings by time $t$ is given as
\begin{align} \label{un_t}
u_{\text{F}}^{\left( n\right) }\left( t,x\right) =& \left(\frac{-i \delta}{\veps} \right)^{n} \frac{1}{\left( 2\pi \varepsilon
\right) ^{3m/2}}\int_{\mathbb{R}^{2m}} dq_{0}dp_{0}
 \int_{0<t_{1}<\cdots<t_{n}<t}dT_{n:1} A_{t} \exp \left( \frac{i}{\varepsilon }%
\Theta_{t}  \right),
 \end{align}
which contains an integral with respect to all possible ordered
hopping times $T_{n:1}:=\left( t_{n},\cdots ,t_{1}\right)$. Here 
\begin{equation}\label{theta} 
  \Theta_{t} = S_{t} +P_{t} \cdot \left(
    x-Q_{t} \right) +\frac{i}{2}\lvert x-Q_{t} \rvert ^{2}.
\end{equation}
and the FGA coefficients $(Q_t, P_t, S_t, A_t)$ corresponds to the position, momentum, phase, and complexed amplitude are governed by  system of ODEs to be described below.

It is shown in \cite{fang-lu} that the series expansion \eqref{eq:u_series} converges absolutely for any
fixed $\veps>0$ and $\delta>0$, since we only integrate over ordered time sequence and thus the summand has a combinatorial factor in the denominator as the number of hops increases, which leads to convergence.

  Each trajectory evolves in the extended phase space and is thus
characterized by position $Q_t$, momentum $P_t$, and surface index
$l_t$ to keep track of the current potential energy surface. The
equation of motion of $(Q, P)$ is given by
\begin{align} 
  \frac{\rd}{\rd t}Q & =P,  \label{ode1} \\
  \frac{\rd}{\rd t}P & =-\nabla V_{l_{t}}(Q) \notag,
\end{align}
where the force depends on the current surface $l_t$. The surface
index $l_{t}$ follows a Markov jump process on the state space
$\left\{0,1\right\}$, where infinitesimal transition probability is given
by%
\begin{equation}
P\left( l_{t+ \tau}=k \ \lvert \ l_{t}=j\right) =\delta _{jk}+\lambda
_{jk} \tau+o\left(\tau\right) ,
\end{equation}%
with transition rate
\begin{equation}
\lambda  =%
\begin{pmatrix}
\lambda_{00}   & \lambda _{01}  \\ 
\lambda_{10}   & \lambda _{11}  
\end{pmatrix}%
=%
\begin{pmatrix}
  -\delta/\eps & \delta/\eps   \\
  \delta/\eps & - \delta/\eps 
\end{pmatrix}%
.  \label{gen}
\end{equation}%
Therefore, $(Q_{t},P_{t},l_{t})$ can be viewed as a
Markov switching process whose probability distribution
$F_{t}\left(q,p,l\right) $ follows the following Kolmogorov forward
equation
\begin{align} \label{eq:kol}
\partial _{t}F_{t}\left(q,p,0\right) + p \cdot \nabla _{q}F_{t}\left(q,p,0\right) -\nabla _{p}V_{0}(q) \cdot \nabla _{p}F_{t}\left(q,p,0\right)=\frac{\delta}{\eps}\left( F_{t}\left(q,p,1\right) - F_{t}\left(q,p,0\right) \right) , \\
\partial _{t}F_{t}\left(q,p,1\right) + p \cdot \nabla _{q}F_{t}\left(q,p,1\right) -\nabla _{p}V_{1}(q) \cdot \nabla _{p}F_{t}\left(q,p,1\right)=\frac{\delta}{\eps}\left( F_{t}\left(q,p,0\right) - F_{t}\left(q,p,1\right) \right),
\end{align}%
where the last two terms of the left hand side is induced by the
Hamiltonian flow of $(Q_{t},P_{t})$, and the right hand side stands
for the contribution of the Markov jump process. Note that $l_t$ is a
piecewise constant and almost surely each trajectory has a finite number
of jumps. Consider a realization of a trajectory with $n$ jumps, we
will denote the discontinuity set of $l_{t}$ as
$\left\{ t_{j}\right\} _{j=1}^{n}$, which is the ordered set of time
of jumps of the trajectory.

Along each trajectory, we also solve for the classical action $S_t$
and weighting factor $A_t$ given by
\begin{align}  \label{ode2}
  \frac{\rd}{\rd t}S & =\frac{1}{2}\lvert P\rvert ^{2}-V_{l_{t}}\left( Q\right) ,   \\
  \frac{\rd}{\rd t}A & =
                  \frac{1}{2}A\, \tr\left(
                  Z^{-1} \left( \partial _{q}P-i\partial _{p} P-i \partial _{q}  Q\nabla _{Q}^{2}{V}_{\ell_t}\left(
                  Q\right) + \partial _{p} Q\nabla _{Q}^{2}{V}_{\ell_t}\left(
                  Q\right)\right) \right),  \label{ode_A}
\end{align}
where $Z:= \partial_q Q + \partial_p P  + i(\partial_q P - \partial_p Q)$.
Throughout our paper, we assume that $P,Q$ and $p,q$ are row vectors, and for $Q = (Q_1, \cdots, Q_d)$ and $q = (q_1, \cdots, q_d)$, the matrix $\partial_q Q$ is defined as
\begin{equation} \label{eq:dqQ}
\partial_q Q = \begin{pmatrix}
\partial_{q_1} Q_1 & \partial_{q_1} Q_2 & \cdots & \partial_{q_1} Q_d \\
\partial_{q_2} Q_1 & \partial_{q_2} Q_2 & \cdots & \partial_{q_2} Q_d \\
\vdots & \vdots & \ddots & \vdots \\
\partial_{q_d} Q_1 & \partial_{q_d} Q_2 & \cdots & \partial_{q_d} Q_d
\end{pmatrix}.
\end{equation}
The other matrices $\partial_p Q$, $\partial_q P$ and $\partial_q Q$ are defined in a similar fashion.
For the trajectory with initial motion
$\left( Q_0, P_0 \right)= \left( q_0, p_0 \right)$, the initial action
 and  amplitude are assigned as
\begin{equation}  \label{initial_A}
S_0 = 0 \qquad \text{and} \qquad 
A_{0} %
=2^{m/2}\int_{\mathbb{R} ^{m}}u_{0}\left(0,
y\right) e^{\frac{i}{\varepsilon }\left( -p_0\left( y-q_0\right) +\frac{i}{2}%
\lvert y-q_0\rvert ^{2}\right) }dy, 
\end{equation}
where the integral for $A_0$ is either analytically given or
numerically estimated.

Finally, the wave function is reconstructed via a trajectory average, that is, an expectation over all of the trajectories,
\begin{equation}  \label{avg1}
u_\text{F}\left( t,x\right) =C_{\mathcal{N}} \mathbb{E}\left[ 
 \left(-i  \right)^n 
 \frac{A_{t}}{\lvert
A_{0} \rvert } \exp \left( \frac{i}{%
\varepsilon }\Theta _{t} \right) 
\exp \left( \frac{\delta}{\eps} t\right) 
\begin{pmatrix}
\mathbb{I}_{n\ \text{is even}} \\ 
\mathbb{I}_{n\ \text{is odd}}%
\end{pmatrix}
\right] ,
\end{equation}
where the normalization constant reads 
\begin{equation}  \label{cn}
C_{\mathcal{N}}=\frac{1}{\left( 2\pi \varepsilon \right) ^{3m/2}}\int_{%
\mathbb{R}^{2m}}\lvert A_{0}\left(q_{0},p_{0}\right)
\rvert dq_{0}dp_{0},
\end{equation}
and we have used $\mathbb{I}$ to denote the indicator function and $\Theta_t$ as defined in \eqref{theta}.

While the expression \eqref{avg1} appears complicated, all the terms
involved in the expectation are obtained along the trajectories. In
the algorithm, the expectation over path ensemble is estimated by
Monte Carlo algorithm sampling in the path space.

\section{The Weak-coupling Regime and Marcus' Rate} \label{sec:weak}
In this section, we consider the non-adiabatic (weak-coupling) regime, namely, $\delta \to 0$ with $\eps$ fixed. In this regime, the celebrated Marcus formula predicts the golden-rule prefactor to be quadratically dependent on the off-diagonal coupling term $\delta$. To be specific, given the initial wavepacket concentrates on only the $0$-th surface and the initial population $\int_{\mathbb{R}^d} \lvert u_\text{in}(x) \rvert^2 dx = 1$, the transition rate $k$ has the following scaling
$$k := \int_{\mathbb{R}^d} \lvert u_1(t,x) \rvert^2 \,d x \sim \frac{\delta^2}{\eps},$$
and the celebrated high-temperature Marcus expression \cite{nitzan2013chemical} reads
$$k = \frac{2\pi \delta^2}{\eps \sqrt{4\pi E_r k_B T}} \exp\left(-\frac{(E_r - \Delta E)^2}{4E_r k_B T} \right),$$
where $E_r$ is the reorganization energy, $\Delta E$ the energy difference between the two diabatic surfaces, $k_B$ the Boltzmann constant and $T$ is the temperature. A exact comparison of the terms require a careful non-dimensionalization. Here we only aim to show the correct scaling in terms of $\delta$ and $\eps$ of the golden-rule prefactor, namely, $O(\delta^2 /\eps)$.
\subsection{A revisit of the deterministic FGA ansatz}
Consider the diabatic FGA ansatz as given in \eqref{eq:u_series} and \eqref{un_t}.
In the weak-coupling regime, it suffices to show that 
$$\int \lvert u_\text{F}^{(1)} (t,x) \rvert^2 dx \sim O\left(\delta^2/\eps\right).$$
To make the notation clearer, we define the coherent state 
$$g^\eps[q,p](x) = (\pi \eps)^{-d/4} \exp\left(-\frac{1}{2\eps}|x-q|^2 + \frac{i}{\eps}p\cdot (x-q)\right).$$
Denote the phase space element $z = (q,p)$ for short, and it follows from a straightforward calculation that the inner product of two coherent states reads
\begin{equation} \label{inner_coherent}
\langle\, g^\eps[z_1], g^\eps[z_2] \rangle = \exp\left(-\frac{1}{4\eps}|z_1 - z_2|^2 + \frac{i}{2\eps}(p_1+p_2)(q_1-q_2)\right).
\end{equation}
Using this notation, the FGA ansatz with no hop can be written as%
\begin{equation*}
u_\text{F}^{(0)}(t,x) = (2\pi \eps)^{-d} \int_{\mathbb{R}^{2d}} A(t,z) e^{i S(t,z)/\eps} g^\eps[Q(t,z),P(t,z)](x) \langle\, g^\eps[z], u_\text{in}\rangle \,dz,
\end{equation*} 
where $(Q,P)$ follows the classical Hamiltonian flow map \eqref{ode1} and the classical action $S$ and the prefactor $a$ follow the ODEs \eqref{ode2} and \eqref{ode_A}. $u_\text{F}^{(1)}$ can be written in a similar fashion, and we arrive at 
\begin{multline}\label{rate_target0}
\int \lvert u_\text{F}^{(1)} (t,x) \rvert^2 dx = \frac{1}{(2\pi\eps)^{2d}} \frac{\delta^2}{\eps^2}\int_0^t d t_1 \int_0^t d t_2 \int_{\mathbb{R}^{2d}} d z_1 \int_{\mathbb{R}^{2d}} d z_2 
\overline{A(t, t_1, z_1)} A(t,t_2,z_2) 
\\ 
\times e^{\frac{i}{\eps}\left(-S(t,t_1,z_1) +S(t,t_2,z_2) \right)} 
\langle g^\eps\left[Q(t,t_1,z_1),P(t,t_1,z_1)\right], g^\eps\left[Q(t,t_2,z_2),P(t,t_2,z_2)\right]\rangle
\overline{\langle g^\eps[z_1], u_\text{in}\rangle} \langle g^\eps[z_2], u_\text{in}\rangle.
\end{multline}
Here with slight abuse of notations, we use $(t_1, z_1)$ to denote the dummy variable in the integral of the jumping time and the phase space for $\overline{u_\text{F}^{(1)}}$, and $(t_2, z_2)$ to denote the dummy variable in the integral for $u_\text{F}^{(1)}$. Both $t_1$ and $t_2$ correspond to the time when the first hop happens but for trajectories representing $\overline{u_\text{F}^{(1)}}$ and $u_\text{F}^{(1)}$, respectively.

\subsection{The Stationary Phase Analysis}
Consider the initial datum as
$$u_\text{in}(x) = g^\eps[z_0](x),$$
then \eqref{rate_target0} can be written as
\begin{equation} \label{rate_target}
\int \lvert u_\text{F}^{(1)} (t,x) \rvert^2 dx = \frac{1}{(2\pi\eps)^{2d}} \frac{\delta^2}{\eps^2}\int_0^t d t_1 \int_0^t d t_2 \int_{\mathbb{R}^{2d}} d z_1 \int_{\mathbb{R}^{2d}} d z_2 \overline{A(t, t_1, z_1)} A(t,t_2,z_2)  e^{\frac{i}{\eps} \Phi(t,t_1,t_2,z_1,z_2)},
\end{equation}
where the phase function is 
\begin{align} \label{eq:phase}
\Phi = & -S(t,t_1,z_1) + S(t,t_2,z_2) 
- \frac{1}{2} (p_1 + p_0)\cdot (q_1-q_0) + \frac{1}{2} (p_2 + p_0)\cdot (q_2-q_0) \\
&+ \frac{1}{2}\left(P(t,t_1,z_1)+P(t,t_2,z_2)\right)\cdot \left(Q(t,t_1,z_1)-Q(t,t_2,z_2)\right) \notag \\
& +\frac{i}{4}|Q(t,t_1,z_1)-Q(t,t_2,z_2)|^2 +\frac{i}{4}|P(t,t_1,z_1)-P(t,t_2,z_2)|^2  \notag\\
& +\frac{i}{4}|z_1-z_0|^2 +\frac{i}{4}|z_2-z_0|^2 \notag.
\end{align} 
Note that \eqref{rate_target} is an oscillatory integral so that the stationary phase asymptotics (see, e.g., \cite{zworski_book}) can be applied, which essentially states that the oscillations of the oscillatory integral act constructively near the non-degenerate stationary point of the phase and destructively elsewhere. To apply the stationary phase asymptotics, we first identify the stationary point of the phase function $\Phi$, defined as the point $(t_1, z_1, t_2, z_2)$ at which the first derivatives of $\Phi$ vanish, i.e.,
\[
\partial_{(t_1, z_1, t_2, z_2)} \Phi (t_1, z_1, t_2, z_2) = 0.
\]
The non-degeneracy of the stationary point will be discussed in \cref{sec:hess}.

\begin{proposition}
The stationary point of the phase function $\Phi$ given by \eqref{rate_target} is the point where
\begin{equation} \label{eq:stationary_pt}
z_1 = z_2 = z_0, \quad t_1 = t_2 \quad \text{such that} \quad V_0(Q(t_1,z_0)) = V_1(Q(t_1,z_0)).
\end{equation}
\end{proposition}
\begin{proof}
To find the stationary points, we set
$$\Im \Phi = 0 \quad \text{and} \quad \nabla_{(t_1,t_2,z_1,z_2)}\Re \Phi=0.$$
The first equation implies 
\begin{gather*} 
z_1 = z_2 = z_0, \\  
Q(t,t_1,z_0)=Q(t,t_2,z_0),  \quad P(t,t_1,z_0)=P(t,t_2,z_0).
\end{gather*}
In other words, the canonical transformation associated with $t_1$ and $t_2$ are the same, and hence $t_1 = t_2$. It is easy to check that $\nabla_{(z_1,z_2)}\Re \Phi=0$ under above conditions, where we used the fact of the classical action that 
\begin{displaymath}
\partial_{q_{\ell}} S(t,t_{\ell},z_{\ell}) = -p_{\ell}^T + [\partial_{q_{\ell}} Q(t,t_{\ell},z_{\ell})] P(t,t_{\ell},z_{\ell})^T, \quad \partial_{p_{\ell}} S(t,t_{\ell},z_{\ell}) =  [\partial_{p_{\ell}} Q(t,t_{\ell},z_{\ell})]P(t,t_{\ell},z_{\ell})^T, \quad \ell=1,2.
\end{displaymath}
Here we regard $\partial_{q_\ell} S$ and $\partial_{p_\ell} S$ as column vectors following the style of \eqref{eq:dqQ}.
It suffices to calculate the partial derivatives on time $t_1$ and $t_2$. For the $t_1$ derivative, we have
\begin{align*}
\partial_{t_1} \Re \Phi =& - \partial_{t_1} S(t,t_1,z_1) 
+\frac{1}{2}\partial_{t_1}P(t,t_1,z_1)\cdot \left(Q(t,t_1,z_1)-Q(t,t_2,z_2)\right) \\
&+\frac{1}{2}\left(P(t,t_1,z_1)+P(t,t_2,z_2)\right)\cdot \partial_{t_1}Q(t,t_1,z_1).
\end{align*}
Evaluating at the critical points $z_1 = z_2 = z_0$, one only needs to show that when $V_0(Q(t_1,z_0)) = V_1(Q(t_1,z_0))$, it holds that
$$\partial_{t_1} S(t,t_1,z_0) =P(t,t_1,z_0)\cdot \partial_{t_1}Q(t,t_1,z_0). $$
Hereafter we drop the dependence on $z_0$ for simplicity. It follows from straightforward calculations that
\begin{gather*}
\partial_{t_1} Q(t,t_1) = \int_{t_1}^t \partial_{t_1} P(\tau,t_1) \, d\tau,\\
\partial_{t_1} S(t,t_1) =  - V_0 (Q(t_1)) + V_1(Q(t_1)) 
+\int_{t_1}^t P(\tau,t_1) \cdot \partial_{t_1} P(\tau,t_1) - \nabla V_1(Q(\tau,t_1)) \cdot \partial_{t_1}Q(\tau,t_1) \,d\tau.
\end{gather*}
Note that in fact
\[
\int_{t_1}^t P(\tau,t_1) \cdot \partial_{t_1} P(\tau,t_1) - \nabla V_1(Q(\tau,t_1)) \cdot \partial_{t_1}Q(\tau,t_1) \,d\tau = P(t,t_1)\cdot \partial_{t_1}Q(t,t_1),  
\]
because
\begin{align*}
& \int_{t_1}^t P(\tau,t_1) \cdot \partial_{t_1} P(\tau,t_1) - \nabla V_1(Q(\tau,t_1)) \cdot \partial_{t_1}Q(\tau,t_1) \,d\tau - P(t,t_1)\cdot \partial_{t_1}Q(t,t_1)  \notag \\
=&  
\int_{t_1}^t \left(P(\tau,t_1) - P(t,t_1) \right) \cdot \partial_{t_1} P(\tau,t_1)
- \nabla V_1(Q(\tau,t_1)) \cdot \partial_{t_1}Q(\tau,t_1) \,d\tau \notag \\
=& \int_{t_1}^t \left(\int_\tau^t \nabla V_1(Q(s,t_1))ds \right) \cdot \partial_{t_1} P(\tau,t_1)
- \nabla V_1(Q(\tau,t_1)) \cdot \left(\int_{t_1}^\tau \partial_{t_1} P(s,t_1) ds \right) \,d\tau  =0,
\end{align*}
where switching the order of the integration contributes to the cancellation in the last line. 
Therefore, we have
\begin{equation} \label{eq:S_t1}
    \partial_{t_1} S(t,t_1) =  - V_0 (Q(t_1)) + V_1(Q(t_1)) + P(t,t_1)\cdot \partial_{t_1}Q(t,t_1),
\end{equation}
which is $P(t,t_1)\cdot \partial_{t_1}Q(t,t_1)$ when applying the stationary condition. The partial derivative of $t_2$ follows a similar calculation. This completes the proof. 
\end{proof}

\begin{remark} 
We note two interesting observations for this lemma. For one thing, the trajectories with $z_1 = z_2 = z_0$ contributes the leading order, and hence it justifies the approach to initialize all trajectories using a fixed $(q_0, p_0)$ rather than using initial sampling, which is typically done in the chemistry literature~\cite{tully1971,tully1990,tully1998mixed,tully2012perspective,landry2011communication,conical_review,WangAkimovPrezhdo2016,SubotnikJainEtAl2016}.  Moreover, the leading contribution comes from hops happening when $ V_0 (Q(t_1)) = V_1(Q(t_1))$, that is, when the adiabatic surfaces achieve minimum gap. So in this scenario, the diabatic surface hopping algorithm that contains multiple jumps can be approximated via a single-jump protocol that can be viewed as a single-jump trajectory surface hopping algorithm~\cite{Lasser2005,lasser2007,fer_las,BelyaevLasserTrigila2014} to some extent. 
\end{remark}

For ease of notation, we denote the Hessian of $\Phi$ with respect to $(t_1, z_1, t_2, z_2)$ at the stationary point \eqref{eq:stationary_pt} as $\hess $.
As long as $\hess$ is non-degenerate,  the principle of stationary phase can be applied and
the transition rate \eqref{rate_target} becomes
\begin{equation*}
\int \lvert u_\text{F}^{(1)} (t,x) \rvert^2 dx = \frac{2\pi \delta^2}{\eps} 
\lvert a(t, t_1, z_0)\rvert^2 (\det \hess )^{-\frac{1}{2}} + O(\delta^2),
\end{equation*}
where $t_1$ is defined in \cref{eq:stationary_pt}. Indeed the scaling factor with respect to $\delta$ and $\epsilon$ is $\delta^2/\epsilon$ that agrees with the Marcus golden-rate scaling. For the transition rate, we have
\begin{equation} \label{eq:all_terms}
    \int_{\mathbb{R}^d} \lvert u_1(t,x) \rvert^2 dx = \sum_{n,m \text{ odd}}^{\infty}  \int_{\mathbb{R}^d} \overline{u^{(n)}_\text{F}(t,x)} u^{(m)}_\text{F}(t,x) \, dx,
\end{equation}
where the term associated with $(n,m)$ contains $n+m$ layers of time integrals, whose leading order read $O \left( (\frac{\delta^2}{\eps})^\frac{n+m}{2}\right)$. Hence the main contribution of the transition rate comes from $\int \lvert u_\text{F}^{(1)} (t,x) \rvert^2 dx$. Thus, it suffices to establish that the Hessian is indeed non-degenerate, as in the next section.

\subsection{Non-degeneracy of the Hessian matrix} \label{sec:hess}
Our main result of this section is to identify a set of generic non-degeneracy conditions of the Hessian matrix at the stationary point:
 \begin{hyp} \label{assump1}
        $\nabla V_0(x) \neq \nabla V_1(x)$ at the points when $V_0(x) = V_1(x)$;
 \end{hyp}
 \begin{hyp} \label{assump2}
        For the time of interest $t > t_1$ with $t_1$ defined according to \eqref{eq:stationary_pt}, 
        $$\partial_{t_1}Q(t, t_1, z_0) \neq 0 .$$ 
 \end{hyp}
\noindent 
For the dynamics under general PESs satisfying both conditions, the diabatic surface hopping algorthm gives the correct scaling factor of the transition rate that matches the Marcus scaling. 
We remark that \cref{assump1} is a generic condition depending only on the two diabatic potential energy surfaces. It states that when the two energy bands cross each other, their derivatives should not be the same. This is a very physical condition for significant transitions between PESs to happen, which can be viewed as a minimum requirement for such results. It is easy to see that the spin-boson model satisfies the condition. \cref{assump2} states that the change of the position with respect to the jumping time $t_1$ needs to be non-zero at the time of interest $t$ . In other words, it means that when jumping happens at different time, one expects the position of the trajectory at the final time $t$ to be different, which is also quite intuitive. In particular, we show that \cref{assump2} is satisfied for the spin-boson model in \cref{sec:appendix}. 

For the rest of this section, we first establish some important properties of the imaginary and real parts of the Hessian, and then prove the non-degeneracy of the Hessian matrix at the stationary point under the generic non-degeneracy conditions. Here $H_R$ and $H_I$ denote the real and imaginary parts of $\hess$ respectively, and $\hess = H_R + i H_I$.

Let us first consider the imaginary part of the Hessian. 
\begin{lemma}\label{lmm:hess_imag_pos_semidef_gen}
$H_I$ is positive semi-definite. 
\end{lemma}
\begin{proof}
    The imaginary part of the total phase \eqref{eq:phase} is 
    \[
       \Im{\Phi} = \frac{1}{4}|Q(t,t_1,z_1)-Q(t,t_2,z_2)|^2 +\frac{1}{4}|P(t,t_1,z_1)-P(t,t_2,z_2)|^2  +\frac{1}{4}|z_1-z_0|^2 +\frac{1}{4}|z_2-z_0|^2.
    \]
    Consider the $\hess_{t_1, z_1, t_2, z_2} \Im{\Phi}$ at the stationary point as defined in \eqref{eq:stationary_pt}. It can be calculated that 
    \begin{equation} \label{eq:hess_I_gen}
        H_I = \frac{1}{2} \begin{pmatrix}
        0 &  &  & \\
          & I_{2d\times 2d} &  & \\
          &  & 0 &  \\
          &  &  & I_{2d\times 2d}
          \end{pmatrix} + \frac{1}{2} \begin{pmatrix}
          F_1 F_1^\dagger & \\
           & F_2 F_2^\dagger
          \end{pmatrix},
    \end{equation}
    where for $\ell = 1, 2$, $F_\ell$ is defined as
    \begin{equation*}
        F_\ell = \begin{pmatrix}
        \partial_{t_\ell} Q(t,t_\ell,z_\ell) & \partial_{t_\ell} P(t,t_\ell,z_\ell) \\
        \partial_{q_\ell} Q(t,t_\ell,z_\ell)  & \partial_{q_\ell} P(t,t_\ell,z_\ell) \\
        \partial_{p_\ell} Q(t,t_\ell,z_\ell)  & \partial_{p_\ell} P(t,t_\ell,z_\ell)
        \end{pmatrix}.
    \end{equation*}
    Thus, the second term of \eqref{eq:hess_I_gen} is positive semi-definite and so is $H_I$.
\end{proof}

For the real part of the Hessian, we make use of \cref{assump1} and \cref{assump2}. %

\begin{lemma} \label{lmm:hess_real_non_degen}
 When both \cref{assump1} and \cref{assump2} are satisfied, $H_R$ is non-degenerate for $t > t_1$.
\end{lemma}
\begin{proof}
For simplicity, we denote $S_\ell = S(t,t_\ell,z_\ell)$, $Q_\ell = Q(t,t_\ell,z_\ell)$ and $P_\ell = P(t,t_\ell,z_\ell)$ for $\ell = 1, 2$ and according to \eqref{eq:phase} we have 
\begin{equation*}
    \Phi_R: = \Re{\Phi} = - S_1 + S_2 - \frac{1}{2}(p_1 + p_0) \cdot (q_1 - q_0) + \frac{1}{2} (p_2 + p_0) \cdot (q_2 - q_0) + \frac{1}{2}(P_1 + P_2) \cdot (Q_1 - Q_2).
\end{equation*}
Taking the first partial derivatives with respect to $q_1$, we have
\begin{align}\label{eq:partial_q1_HR}
    \partial_{q_1} \Phi_R = & -\partial_{q_1} S_1 + \frac{1}{2} (\partial_{q_1} P_1) (Q_1 - Q_2)^T + \frac{1}{2} (\partial_{q_1} Q_1) (P_1 + P_2)^T - \frac{1}{2}(p_1 + p_0)^T  \notag
    \\
    = & \frac{1}{2} (\partial_{q_1} P_1) (Q_1 - Q_2)^T - \frac{1}{2} (\partial_{q_1} Q_1) (P_1 - P_2)^T + \frac{1}{2}(p_1 - p_0)^T,
\end{align}
where we used the fact that $\partial_{q_1} S = - p_1^T +  (\partial_{q_1} Q_1) P_1^T$. Similarly, its first partial derivative with respect to $p_1$ can be written as
\begin{align} \label{eq:partial_p1_HR}
    \partial_{p_1} \Phi_R = &- \partial_{p_1} S_1 + \frac{1}{2} (\partial_{p_1}P_1) (Q_1 - Q_2)^T + \frac{1}{2} (\partial_{p_1} Q_1) (P_1 + P_2)^T - \frac{1}{2}(q_1 - q_0)^T \notag
    \\
    = & \frac{1}{2} (\partial_{p_1}P_1) (Q_1 - Q_2)^T - \frac{1}{2} (\partial_{p_1} Q_1) (P_1 - P_2)^T - \frac{1}{2}(q_1 - q_0)^T,
\end{align}
thanks to the fact that $\partial_{q_1} S_1 = (\partial_{q_1} Q_1) P_1^T$. We now consider the second partial derivatives of $\Phi_R$ evaluated at the stationary point. Note that at the stationary point $P_1 = P_2$ and $Q_1 = Q_2$, resulting that
\begin{align*}
    \partial_{q_1}^2 \Phi_R = &  \frac{1}{2} (\partial_{q_1} P_1) (\partial_{q_1} Q_1)^T - \frac{1}{2} (\partial_{q_1} Q_1) (\partial_{q_1} P_1)^T,
    \\
    \partial_{q_1} \partial_{q_2} \Phi_R = & -\frac{1}{2}(\partial_{q_1} P_1) (\partial_{q_2} Q_2)^T + \frac{1}{2}(\partial_{q_2} Q_2) (\partial_{q_1} P_1)^T = -\partial_{q_1}^2 \Phi_R.
\end{align*}
Using the fact that $\partial_{q_1}^2 \Phi_R$ is a symmetric matrix, we see that
\[
   \partial_{q_1}^2 \Phi_R = \partial_{q_1} \partial_{q_2} \Phi_R = 0.
\]

A key observation is that for $\ell = 1, 2$,
\begin{equation} \label{eq:QqPp=1}
    (\partial_{q_\ell} Q_\ell) (\partial_{p_\ell} P_\ell)^T - (\partial_{q_\ell} P_\ell) (\partial_{p_\ell} Q_\ell)^T = I.
\end{equation}
This can be derived by considering the derivative of the left-hand-side of \eqref{eq:QqPp=1} with respect to $t$, which is zero thanks to 
\begin{equation} \label{eq:ode_jacobian}
\frac{d}{dt} \begin{pmatrix}\partial_q Q & \partial_q P \\ \partial_p Q & \partial_p P \end{pmatrix} 
=   \begin{pmatrix}\partial_q Q & \partial_q P \\ \partial_p Q & \partial_p P \end{pmatrix} \begin{pmatrix} 0& -\nabla_Q^2 {V}_{\ell_t}(Q)  \\ I & 0\end{pmatrix}.
\end{equation}
In light of \eqref{eq:QqPp=1}, we obtain the partial derivatives involving both $q$ and $p$ at the stationary point
\begin{align*}
    \partial_{q_1}\partial_{p_1} \Phi_R = &  \frac{1}{2} (\partial_{q_1} P_1) (\partial_{p_1} Q_1)^T - \frac{1}{2} (\partial_{q_1} Q_1) (\partial_{p_1} P_1)^T + \frac{1}{2}I = -\frac{1}{2}I + \frac{1}{2}I = 0, 
    \\
    \partial_{q_1}\partial_{p_2} \Phi_R = &  -\frac{1}{2} (\partial_{q_1} P_1) (\partial_{p_2} Q_2)^T + \frac{1}{2} (\partial_{q_1} Q_1) (\partial_{p_2} P_2)^T = \frac{1}{2}I,
\end{align*}
where in the second line we used the stationary condition that $\partial_{p_2} Q_2 = \partial_{p_1} Q_1$ and $\partial_{p_2} P_2 = \partial_{p_1} P_1$. The other second-order partial derivatives of $\Phi_R$ with respect to $(z_1, z_2)$ follows a similar calculation and we are left with the evaluation of the second-order partial derivatives involving $t_1$ or $t_2$. Because
\begin{align*}
    \partial_{t_1} \Phi_R  =  - \partial_{t_1} S_1 + \frac{1}{2} \partial_{t_1} P_1 \cdot (Q_1 - Q_2) + \frac{1}{2}(P_1 + P_2) \cdot \partial_{t_1} Q_1,
\end{align*}
at the stationary point
\begin{equation*}
    \partial_{t_1}\partial_{t_2} \Phi_R = - \frac{1}{2} \partial_{t_1} P_1 \cdot \partial_{t_2} Q_2 + \frac{1}{2} \partial_{t_2} P_2 \cdot \partial_{t_1} Q_1 = 0.
\end{equation*}
Thanks to \eqref{eq:partial_p1_HR}, \eqref{eq:partial_q1_HR} and the stationary condition \eqref{eq:stationary_pt}, we have 
\begin{align*}
    & \partial^2_{t_1} \Phi_R = - \partial^2_{t_2} \Phi_R = - \partial_{t_1}^2 S_1 + P_1 \cdot \partial^2_{t_1}Q_1 +\partial_{t_1} P_1 \cdot \partial_{t_1} Q_1 : = A, 
    \\
    & \partial_{t_1} \partial_{q_1} \Phi_R = \partial_{t_1} \partial_{q_2} \Phi_R = - \partial_{t_2} \partial_{q_1} \Phi_R =  - \partial_{t_2} \partial_{q_2} \Phi_R : = B,
    \\
    & \partial_{t_1} \partial_{p_1} \Phi_R  = \partial_{t_1} \partial_{p_2} \Phi_R = - \partial_{t_2} \partial_{p_1} \Phi_R =  - \partial_{t_2} \partial_{p_2} \Phi_R : = C.
\end{align*}
We are ready to discuss the determinant of $H_R$. Let $$M = \begin{pmatrix} A & B^T & C^T \\ B & 0 & 0 \\ C & 0 & 0 \end{pmatrix}, \qquad N = \begin{pmatrix} 0 & B^T & C^T \\ -B & 0 & \frac{1}{2} I \\ -C & -\frac{1}{2} I & 0 \end{pmatrix}.$$ Then 
\begin{align*}
    \det H_R = & \det \begin{pmatrix} M & N \\ -N & -M \end{pmatrix} = \det \left[ \begin{pmatrix} I & 0 \\ I & I \end{pmatrix} \begin{pmatrix} M & N \\ -N & -M \end{pmatrix} \begin{pmatrix} I & 0 \\ I & I \end{pmatrix} \right] 
    \\
    = & -\det(M+N) \det(M-N) = -\frac{1}{2^{4d}} A^2.
\end{align*}
Taking advantage of \eqref{eq:S_t1}, a further calculation reveals that
\begin{align*}
A = & \partial_{t_1} \left( V_0 (Q( t_1)) - V_1(Q(t_1)) - P_1\cdot \partial_{t_1}Q_1\right) + P_1 \cdot \partial^2_{t_1}Q_1 + \partial_{t_1} P_1 \cdot \partial_{t_1} Q_1 
\\
= &  \left(\nabla V_0 (Q(t_1)) - \nabla V_1(Q(t_1))\right) \cdot \partial_{t_1}Q_1,
\end{align*}
where $Q(t_1) = Q(t_1, z_0)$ is the position at the jumping time $t_1$ that satisfies the stationary condition \eqref{eq:stationary_pt}.
Thanks to \cref{assump1} and \cref{assump2}, $\nabla V_0 (Q(t_1)) - \nabla V_1(Q(t_1) \neq 0$ and $\partial_{t_1}Q_1 = \partial_{t_1}Q(t,t_1,z_0) \neq 0 $.
\end{proof}

We are now ready to establish the non-degeneracy of the Hessian matrix. 
\begin{proposition}\label{lmm:hess_non_degen}
    If both \cref{assump1} and \cref{assump2} are satisfied, $\det(\hess) \neq 0 $ for $t>t_1$.
\end{proposition}
\begin{proof}
    We denote the real and imaginary parts of $\hess$ as $H_R$ and $H_I$ that are both real-valued, so that $\hess = H_R + i H_I$. Thanks to \cref{lmm:hess_imag_pos_semidef_gen}, $H_I$ is positive semi-definite and hence there exists some real-valued square matrix $J$ such that 
    \[
    H_I = J J^T.
    \]
    Together with \cref{lmm:hess_real_non_degen}, the determinant of $\hess$ can be written as
    \[
    \det(H_R + i J J^T) = \det(I + i J J^T H_R^{-1}) \det{H_R} = \det(I + i J^T H_R^{-1} J) \det{H_R},
    \]
    where Sylvester's determinant theorem ($\det(I + AB) = \det(I + BA)$, \cite[p. 271]{Pozrikidis2014}) has been applied. Note that $M: = i J^T H_R^{-1} J$ is skew-Hermitian, because
    \[
    M^\dagger = -i J^T (H_R^{-1})^T J = -M,
    \]
    where one uses the fact that $H_R$ is real symmetric. Furthermore $\det{H_R} \neq 0$ thanks to \cref{lmm:hess_real_non_degen}. Therefore, it suffices to show that 
    \[
    \det(I + i J^T H_R^{-1} J) =  \det(I + M) \neq 0.
    \]
    Suppose on the contrary that there exists some vector $v \neq 0$ such that 
    \[
    (I + M) v = 0.
    \]
    Multiplying by $v^\dagger$ on both sides yields that
    \[
    0 = v^\dagger (I + M)  v = \norm{v}^2 + v^\dagger M v = \norm{v}^2,
    \]
    where we used the fact that $v^\dagger M v = 0$ for skew-Hermitian matrix $M$. This is a contradiction, which completes the proof of the desired result.
\end{proof}

\medskip
\section{The Strong-coupling Regime and the Ehrenfest Dynamics} \label{sec:large}
In this section, we study the limiting behavior of the diabatic surface hopping algorithm in the strong-coupling regime $\delta \to \infty$ with $\eps$ fixed. Note that this corresponding to the adiabatic regime, as the gap between the two adiabatic PESs approaches to infinity. In contrast to the limiting behavior of trajectory-based nonadiabatic algorithm in the adiabatic representation, where the dynamics approximately stay on the same adiabatic PESs as it starts in thanks to the Born-Oppenheimer approximation, the surface hopping algorithm in the diabatic picture suffers from a huge jumping rate (large coupling constant), which refrains one from adopting perturbation approaches. Unlike the weak-coupling regime, all terms in the series \eqref{eq:u_series} need to be taken into account. However, the analysis is drastically simplified when considering the stochastic representation \eqref{avg1} instead, which essentially encodes all terms in the series.

Heuristically, one can imagine as the trajectories hop between the two diabatic PESs constantly, the effective diabatic PES may approach to the mean-field (or Ehrenfest) dynamics $\frac{1}{2}(V_0 +V_1)$.  Hence, replacing the propagator in the Dyson series with the mean dynamics gives
\begin{equation} \label{eq:heuristic_ehrenfest}
\begin{pmatrix}
\cos(\frac{\delta}{\eps}t) \bar{\mathcal{U}}(t,0) u_0^\text{in} \\
-i\sin(\frac{\delta}{\eps}t) \bar{\mathcal{U}}(t,0) u_0^\text{in},
\end{pmatrix},
\end{equation}
where $ \bar{\mathcal{U}}$ denotes the propagator associated with the mean potential $\frac{1}{2}(V_0 +V_1)$.

The focus of this section is to justify this intuition for general diabatic potentials. In particular, we show that the mean-field dynamics can be indeed derived in the limit  $\delta / \epsilon \to \infty$ using a two-scale expansion.

\subsection{Estimator of the wave function}

Let us make transparent of the trajectory expectation in the case $\delta/\eps \to \infty$. The estimator of the wave function is given by \eqref{avg1},
$A_0$ defined by \eqref{initial_A}, and $(Q_t, P_t, S_t)$ depends on the history of hopping time $T_{n:1}$ and the initial datum $z$. %
Note that $A_t$ in fact can be represented in terms of the Jacobian matrix of the underlying flow ~\cite{swart2009mathematical}
\[
A_t = \sqrt{2^{-d} \det\left(\partial_q Q + \partial_p P  + i(\partial_q P - \partial_p Q) \right) },
\]
where we recall that the components of the Jacobian matrix $(\partial_q Q, \partial_q P, \partial_p Q, \partial_p P)$ follow the system of ODEs \eqref{eq:ode_jacobian}.
Though this representation involves taking the square root of a complex variable that can be less appealing in numerical implementations, it is helpful for our theoretical study here. We are now ready for the derivation.

Define $\tilde{z}_t : = (z_t, \ell_t)$, where we recall that $z_t$ is the phase space variable and $\ell_t$ is the surface index. Applying the law of total expectation, \eqref{avg1} becomes
\begin{align*} 
u_\text{F}\left( t,x\right )  &=C_{\mathcal{N}} E_{z} \left\{\mathbb{E}_{\tilde{z}_t}\left[ \left.
 \left(-i  \right)^n  
 \frac{A_{t}}{\lvert
A_{0} \rvert } \exp \left( \frac{i}{%
\varepsilon }\Theta _{t} \right)  \exp \left( \frac{\delta t}{\eps}\right) 
\begin{pmatrix}
\mathbb{I}_{n\ \text{is even}} \\ 
\mathbb{I}_{n\ \text{is odd}}%
\end{pmatrix}
 \right\vert  \tilde{z}_0 = (z, 0) \right] \right\}
 \\
&= C_{\mathcal{N}}  \int dz   \mathbb{E}_{\tilde{z}_t}\left[ \left.
 \left(-i  \right)^n   A_t(z) \exp \left( \frac{i}{%
\varepsilon }\Theta _{t} \right)  \exp \left( \frac{\delta t}{\eps}\right) 
\begin{pmatrix}
\mathbb{I}_{n\ \text{is even}} \\ 
\mathbb{I}_{n\ \text{is odd}}%
\end{pmatrix}
 \right\vert  \tilde{z}_0 = (z, 0) \right] \\
 & = C_{\mathcal{N}}  \int dz    \int dX_t \sum_n G(t, X_t, n)
 \left(-i  \right)^n \alpha(X_t)
 \exp \left( \frac{i}{%
\varepsilon }\Theta _{t} \right) \exp \left( \frac{\delta t}{\eps}\right) 
\begin{pmatrix}
\mathbb{I}_{n\ \text{is even}} \\ 
\mathbb{I}_{n\ \text{is odd}}%
\end{pmatrix}
,
\end{align*}
where $G(t,X,n)$ is the joint distribution of $$X =(Q, P, S,\partial_q Q, \partial_q P, \partial_p Q, \partial_p P)$$
and $n$ at time $t$, and $\alpha$ is defined as \[
\alpha(X) = \sqrt{2^{-d} \det\left(\partial_q Q + \partial_p P  + i(\partial_q P - \partial_p Q) \right)}.
\] We write $F_i(t,X) = \sum_{n \in \mathcal{I}_i} G(t,X,n)$ for $i = 0, 1, 2, 3$ and $\mathcal{I}_i = \{n : n \bmod 4 = i\}$. We have
\begin{equation*} 
u_\text{F}\left( t,x\right)
  = C_{\mathcal{N}}  \int dz 
 \int dX_t  \alpha(X_t) \exp \left( \frac{i}{%
\varepsilon }\Theta _{t} \right) \exp \left( \frac{\delta t}{\eps}\right) 
\begin{pmatrix}
F_0 (t, X_t) - F_2(t, X_t) \\ 
(-i)\left(F_1 (t, X_t) - F_3(t, X_t)\right)
\end{pmatrix}
.
\end{equation*}
Note that the results depend only on the differences
\[
\Delta F_0 : = F_0 - F_2, \quad \Delta F_1 : = F_1-F_3
\]
that satisfy the following equations
\begin{equation} \label{eq:DF}
\begin{aligned}
    {\partial \over \partial t} \Delta F_0(t,X) + \nabla_X \cdot(W_0(X) \Delta F_0(t,X)) &= \frac{\delta}{\eps}\left( \Delta F_1(t,X) - \Delta F_0(t,X)\right), \\
    {\partial \over \partial t} \Delta F_1(t,X) + \nabla_X \cdot(W_1(X)\Delta F_1(t,X)) &= \frac{\delta}{\eps}\left(- \Delta F_0(t,X) - \Delta F_1(t,X)\right).
\end{aligned}
\end{equation}
Here $W_\ell$ (for $\ell = 0,1$) is given by
\[
W_\ell(Q,P,S) = \left( P, \  -\nabla V_\ell(Q),  \ \frac{1}{2}|P|^2 -  V_\ell(Q), \ \partial_q P, \  - \partial_q Q \, \nabla _{Q}^{2}V_{\ell} (Q), \ \partial_p P, \ -\partial_p Q \,\nabla _{Q}^{2}V_{\ell}(Q) \right),
\]
and the intial data are given as 
\begin{align} \label{eq:initial_Delta_F0_F1}
    \Delta F_0(0,Q,P,S,\partial_q Q, \partial_q P, \partial_p Q, \partial_p P) = & \delta\left(Q=q, P=p, S=0,\partial_q Q = I, \partial_q P = 0, \partial_p Q = 0, \partial_p P = I\right), 
    \\
    \Delta F_1(0,Q,P,S,\partial_q Q, \partial_q P, \partial_p Q, \partial_p P)= & 0. \notag
\end{align}
\subsection{A two-scale expansion}
Now we are going to find asymptotic solutions to \eqref{eq:DF}. Define
\begin{equation*}
f_0(t,X) = \exp \left( \frac{\delta t}{\eps} \right) \Delta F_0(t,X), \qquad
f_1(t,X) = \exp \left( \frac{\delta t}{\eps} \right) \Delta F_1(t,X).
\end{equation*}
By \eqref{eq:DF}, the functions $f_0$ and $f_1$ satisfy the equations
\begin{equation} \label{eq:f}
\begin{aligned}
{\partial \over \partial t} f_0(t,X) + \nabla_X \cdot(W_0(X) f_0(t,X)) &= \frac{\delta}{\eps} f_1(t,X), \\
{\partial \over \partial t} f_1(t,X) + \nabla_X \cdot(W_1(X) f_1(t,X)) &= -\frac{\delta}{\eps} f_0(t,X).
\end{aligned}
\end{equation}
Due to the smallness of $\eps/\delta$, the advection term and the right-hand side of the above system indicate, respectively, a slow and fast variation with time. Therefore we introduce the fast scale $t^\ast = \delta t / \eps$, and assume
\begin{equation*}
f_0(t,X) = f_0^{(0)}(t,t^\ast,X) + {\eps \over \delta} f_0^{(1)}(t,t^\ast,X) + \cdots, \qquad
f_1(t,X) = f_1^{(0)}(t,t^\ast,X) + {\eps \over \delta} f_1^{(1)}(t,t^\ast,X) + \cdots.
\end{equation*}
By inserting the above ansatz into \eqref{eq:f} and matching the $O(\delta/\eps)$ terms, we obtain
\begin{equation*}
\begin{aligned}
{\partial \over \partial t^\ast} f_0^{(0)}(t,t^\ast,X) &= f_1^{(0)}(t,t^\ast,X), \\
{\partial \over \partial t^\ast} f_1^{(0)}(t,t^\ast,X) &= -f_0^{(0)}(t,t^\ast,X).
\end{aligned}
\end{equation*}
The solution is
\begin{equation} \label{eq:zeroth}
\begin{pmatrix} f_0^{(0)}(t,t^\ast,X) \\ f_1^{(0)}(t,t^\ast,X) \end{pmatrix}
  = R(t^\ast) \begin{pmatrix} a(t,X) \\ b(t,X) \end{pmatrix},
\end{equation}
where the matrix $R(t)$ is
\begin{equation*}
R(t^\ast) = \begin{pmatrix} \cos t^\ast & \sin t^\ast \\ -\sin t^\ast & \cos t^\ast \end{pmatrix}
\end{equation*}
The $O(1)$ equations are
\begin{equation*}
R(t^\ast) {\partial \over \partial t} \begin{pmatrix} a(t,X) \\ b(t,X) \end{pmatrix}
  + {\partial \over \partial t^\ast}
    \begin{pmatrix} f_0^{(1)}(t,t^\ast,X) \\ f_1^{(1)}(t,t^\ast,X) \end{pmatrix}
  + \nabla_X \cdot \left[ \overline{W}(X) R(t^\ast)
    \begin{pmatrix} a(t,X) \\ b(t,X) \end{pmatrix} \right]
= \begin{pmatrix} f_1^{(1)}(t,t^\ast,X) \\ -f_0^{(1)}(t,t^\ast,X) \end{pmatrix},
\end{equation*}
where $\overline{W}(X) = \mathrm{diag} \{ W_0(X), W_1(X) \}$. The above equation can be reformulated as
\begin{equation} \label{eq:ab}
{\partial \over \partial t} \begin{pmatrix} a(t,X) \\ b(t,X) \end{pmatrix}
  + {\partial \over \partial t^\ast} \left[ \big( R(t^\ast) \big)^{-1}
    \begin{pmatrix} f_0^{(1)}(t,t^\ast,X) \\ f_1^{(1)}(t,t^\ast,X) \end{pmatrix} \right]
  + \nabla_X \cdot \left[ \big( R(t^\ast) \big)^{-1} \overline{W}(X) R(t^\ast)
    \begin{pmatrix} a(t,X) \\ b(t,X) \end{pmatrix} \right] = 0.
\end{equation}
To determine $a(t,X)$ and $b(t,X)$, we further assume that the solution has a periodic microstructure, so that $f_0^{(k)}(t,t^\ast,X)$ and $f_1^{(k)}(t,t^\ast,X)$ are periodic with respect to $t^\ast$. The period $T = 2\pi$ can be read off from the zeroth-order solution \eqref{eq:zeroth}. Thus, by integrating \eqref{eq:ab} from $0$ to $T$ with respect to $t^\ast$, we get the equations for $a(t,X)$ and $b(t,X)$:
\begin{equation} \label{eq:a_b}
{\partial \over \partial t} \begin{pmatrix} a(t,X) \\ b(t,X) \end{pmatrix}
  + \nabla_X \cdot \left[
    \frac{1}{2} \begin{pmatrix} W_0(X) + W_1(X) & 0 \\ 0 & W_0(X) + W_1(X) \end{pmatrix}
    \begin{pmatrix} a(t,X) \\ b(t,X) \end{pmatrix} \right] = 0,
\end{equation}
which matches the mean-field Ehrenfest dynamics governed by the average potential $\frac{1}{2}\left(W_0 + W_1 \right)$ as desired.

We now take into account the initial condition \eqref{initial_data} that concentrates only on the first diabatic PES. Because of \eqref{eq:initial_Delta_F0_F1} and \eqref{eq:zeroth},
the initial condition of \eqref{eq:a_b} is given as 
\begin{align*}
    a(0,Q,P,S,\partial_q Q, \partial_q P, \partial_p Q, \partial_p P) = & \delta (Q = q, P= p, S = 0,\partial_q Q = I, \partial_q P = 0, \partial_p Q = 0, \partial_p P = I), \\
    b(0,Q,P,S,\partial_q Q, \partial_q P, \partial_p Q, \partial_p P) = & 0.
\end{align*} 
Hence, $b$ stays zero for all time and $a$ remains the structure of a $\delta$-function
\[
a(t, X) = \delta (t, X = \bar X_t)
\]
where $\bar X_t : = (\bar Q_t, \bar P_t, \bar S_t,\overline{\partial_q Q}_t ,  \overline{\partial_q P}_t, \overline{\partial_p Q}_t, \overline{\partial_p P}_t)$ are the quantities of time $t$ following the average dynamics $(W_0(X) + W_1(X))/2$. In other words, $\bar X_t$ are the position, momentum, action and components of the Jacobian matrix that are governed by the Hamiltonian flow with the potential $(V_0 + V_1)/2$.
The wavefunction constructed by the diabatic algorithm up to leading order then becomes
\begin{align*} 
u_\text{F}\left( t,x\right)
  = &  C_{\mathcal{N}}  \int dz  
\int dX_t \alpha(X_t) \exp \left( \frac{i}{%
\varepsilon }\Theta _{t} \right) 
\begin{pmatrix}
1 & 0 \\ 
0 & -i
\end{pmatrix}
\begin{pmatrix} \cos \left(\frac{\delta t}{ \eps}\right) & \sin \left(\frac{\delta t}{ \eps}\right) \\ -\sin \left(\frac{\delta t}{ \eps}\right) & \cos \left(\frac{\delta t}{ \eps}\right) \end{pmatrix} \begin{pmatrix} a(t,X_t) \\ 0 \end{pmatrix}
\\
= & C_{\mathcal{N}}  \int dz  \alpha (\bar X_t) 
 \exp \left( \frac{i}{%
\varepsilon }\bar S_{t} +\bar P_{t} \left(
    x-\bar Q_{t} \right) +\frac{i}{2}\lvert x-\bar Q_{t} \rvert ^{2} \right) 
\begin{pmatrix} 
\cos \left(\frac{\delta t}{ \eps}\right)   \\ -i\sin \left(\frac{\delta t}{ \eps}\right) 
\end{pmatrix}. 
\end{align*}
Notice that the $z$-integral part is the FGA representation of the scalar (single-band) Schr\"odinger wave function with the PES given by $(V_0 + V_1)/2$. Therefore, we recover \eqref{eq:heuristic_ehrenfest} which justifies the heuristics.

We remark that though we focus on the derivation of the mean dynamics, a derivation for the variance dynamics is also possible using
the quasi steady-state diffusion approximation that is developed for transport processes \cite{Papanicolaou1975} and widely applied to biological applications. While we will not go into the details here.

\section{Conclusion and Discussion} \label{sec:conclusion}
We have studied in details the behavior of the diabatic surface hopping algorithm \cite{fang-lu} in both the weak-coupling (non-adiabatic) regime and the strong-coupling (adiabatic) regime. In the weak-coupling regime, we demonstrate that the algorithm captures the correct scaling of the transition rate that matches the celebrated Marcus theory. In particular, the leading term is the contribution from hopping once, and the transition rate can thus be written as a highly oscillatory integral. The proof relies on the stationary phase analysis and perturbation techniques. We proposed two non-degeneracy conditions of the Hessian matrix, under which the transition rates are shown to capture the correct Marcus scaling for general diabatic PESs. These conditions are naturally satisfied by the spin-boson model. Our result matches the numerical evidence in \cite{fang-lu}. It is an interesting question to investigate the stationary phase asymptotics for the Landau-Zener regime $\delta = \sqrt{\epsilon} \to 0$ \cite{hagedorn1998landau,landau1932theorie,zener1932non}, where the asymptotics need to be carried out for all terms in \cref{eq:all_terms}, and see whether one can prove the diabatic surface hopping algorithm captures the correct Landau-Zener formula. It is also interesting to investigate physical observables besides the transition rate, which may be possible by adopting Wely quantization techniques~\cite{Martinez2002,zworski_book,RobertCombescure2021} that have been exploited to obtain observable bounds for the algorithms for the scalar semiclassical Schr\"odinger equation~\cite{LasserRoblitz2010,LasserLubich2020,FangTres2021}.

In the strong-coupling regime, the perturbation perspective no longer holds and contributions from all number of hoppings need to be considered. Nevertheless, the limiting behavior can be derived by utilizing  the stochastic representation instead. Using a two-scale expansion, we find in the mean dynamics the diabatic surface hopping algorithm behaves as a type of mean-field Ehrenfest dynamics governed by an average potential. The techniques we used here for investigating the two limiting regimes can be of potential use in the analysis for other diabatic algorithms. The limiting dynamics can also aid in the variance reduction and improvements of the existing algorithm, which is an interesting future research direction. 

It is interesting to observe that the non-perturbative (and hence difficult) regime in the adiabatic representation is the perturbative (and easy) regime in the diabatic representation and vice visa. Therefore, our analysis suggest that a combination of both pictures in different regimes can be an appealing approach in handling general non-adiabatic dynamics; this will be left for future research.

\appendix 
\section{Verification of non-degeneracy condition for the spin-boson model} \label{sec:appendix}

For the spin-boson model, the non-degeneracy conditions can be verified via direct calculations. It is immediate that \cref{assump1} is satisfied, and we focus on the verification of \cref{assump2}. Consider the time $t$ after the transition time but still within this oscillating period (before returning to the crossing points again), for example, $w(t - t_1) \in (0, \pi)$.
\begin{lemma}[\Cref{assump2} for spin-boson] For time $t$ such that $0< w(t-t_1) < \pi$, \cref{assump2} is satisfied, where $t_1$ is defined as \eqref{eq:stationary_pt}.
\end{lemma}
\begin{proof}

 The ODEs that $(Q,P)$ satisfies can be written as
 \begin{align*}
     \frac{\rd}{\rd t} \begin{pmatrix}
     Q \\ P
     \end{pmatrix} 
     = \begin{pmatrix}
     0 & 1 \\ -w^2 & 0
     \end{pmatrix}  \begin{pmatrix}
     Q \\ P
     \end{pmatrix} 
     + 
      \begin{pmatrix}
     0 \\ c(t)
     \end{pmatrix} ,
 \end{align*}
 where the vector $c(t) = -c$ when $t \in [0, t_1]$ and $c(t) = c$ when $t \geq t_1$. Denote $\begin{pmatrix}
     0 & 1 \\  -w^2 & 0
     \end{pmatrix}$ as $M_0$. Applying the variation of constant formula, we have
 \begin{align*}
     \begin{pmatrix}
     Q \\ P
     \end{pmatrix}  
     =  e^{t M_0} \begin{pmatrix}
     q_0\\ p_0
     \end{pmatrix} - \int_0^{t_1} e^{ (t -s) M_0} \begin{pmatrix}
     0 \\ c 
     \end{pmatrix} \, ds + \int_{t_1}^{t} e^{ (t -s) M_0} \begin{pmatrix}
     0 \\ c
     \end{pmatrix}\, ds,
 \end{align*}
 which gives the expression of $Q$ via simple calculation  
 \begin{equation*}
     Q(t, t_1, z_0) =\left (q_0  + \frac{c}{w^2} \right)\cos(wt) + \frac{p_0}{w}\sin(wt) + \frac{c}{w^2} - \frac{2c}{w^2}\cos(w(t-t_1)).
 \end{equation*}
 Taking the partial derivative with respect to $t_1$ yields 
 \begin{equation*} \label{eq:spin_boson_partial_t_Q}
    \partial_{t_1} Q(t, t_1, z_0)=  -\frac{2c}{w} \sin(w(t-t_1)) \neq 0,
\end{equation*}
which finishes the proof.
\end{proof}

\section*{Acknowledgments} 
The work of Z.C. was supported by the Academic Research Fund of the Ministry of Education of Singapore under grants R-146-000-291-114 and R-146-000-326-112.
D.F. acknowledges the support by the NSF Quantum Leap Challenge Institute (QLCI) program through grant number OMA-2016245 and by Department of Energy through the Quantum Systems Accelerator program and the Simons Institute Quantum research pod. The work of J.L. is supported in part by National Science Foundation via grant DMS-2012286 and by Department of Energy via grant DE-SC0019449. 

\bibliographystyle{abbrvurl}
\bibliography{sh}

\end{document}